\documentclass[11pt]{amsart}

\usepackage{amssymb}

\def\P{{\mathcal{P}}}
\def\F{{\mathcal{F}}}
\def\T{{\mathcal{T}}}

\def\Nn{{\mathbb{N}}}

\let\notto\nrightarrow

\usepackage[arrow,matrix]{xy}

\makeatletter
\def\aliasbib#1#2{\expandafter\xdef
\csname b@#1\endcsname{\csname b@#2\endcsname}}
\makeatother

\aliasbib{pultr}{Pul:The-right-adjoints}
\aliasbib{tarmul}{Tar:Mul}
\aliasbib{HN:Homomorphisms}{HelNes:GrH}
\aliasbib{es}{EZaSau:Chro}
\aliasbib{Tar:dbl}{Tar:Hedetniemis}

\title{Adjoint functors in graph theory}
\author[J. Foniok]{Jan Foniok}
\author[C. Tardif]{Claude Tardif}

\date{April 8, 2013}

\address{Department of Mathematics and Statistics \\
Queen's University\\
Jeffery Hall\\
48 University Avenue\\
Kingston ON K7L 3N6\\
Canada}
\email{foniok@mast.queensu.ca}

\address{Royal Military College of Canada \\ 
PO Box 17000 Station ``Forces'' \\
Kingston, Ontario\\ 
Canada, K7K 7B4 } 
\email{Claude.Tardif@rmc.ca}
\urladdr{http://www.rmc.ca/academic/math\underline{\ }cs/tardif/}

\thanks{The second author's research is supported 
by grants from NSERC and ARP}


\usepackage{amsthm}
\newtheorem{lemma}{Lemma}[section]
\newtheorem{theorem}[lemma]{Theorem}

\newtheorem{prop}[lemma]{Proposition}
\theoremstyle{definition}
\newtheorem{defi}[lemma]{Definition}
\newtheorem{problem}[lemma]{Problem}
\newtheorem*{exas}{Examples}
\newtheorem*{exa}{Example}

\newfont{\Bb}{msbm10 scaled\magstep1}

\begin{document}


\begin{abstract}
We survey some uses of adjoint functors in graph theory pertaining to
colourings, complexity reductions, multiplicativity, circular
colourings and tree duality. The exposition of these applications through
adjoint functors unifies the presentation to some extent, and also
raises interesting questions. 
\end{abstract}

\date{\today}
\maketitle

\section{Introduction}
We will motivate our subject with known examples from the literature,
which the reader may recognize. Relevant definitions are postponed 
to the following section.

\begin{theorem}[Geller, Stahl] \label{gest}
Let $G$ be a graph. Then the lexicographic
product $G[K_2]$ is $n$-colourable if and only if $G$ 
admits a homomorphism to the Kneser graph $K(n,2)$.
\end{theorem}

\begin{theorem}[El-Zahar, Sauer] \label{elsa}
For two graphs $G$ and $H$, the categorial product $G\times H$
is $n$-colourable if and only if $G$ admits a homomorphism to 
the exponential graph $K_n^{H}$.
\end{theorem}

\begin{theorem}[Hell, Ne\v set\v ril] \label{albe} For a graph $G$, let $G^{1/3}$ be the graph
obtained by replacing each edge of $G$ by a path with three edges.
Then there exists a homomorphism of $G^{1/3}$ to the $5$-cycle $C_5$
if and only if $G$ is $5$-colourable.
\end{theorem}

An $n$-colouring of a graph $G$ is a homomorphism of $G$ to the complete graph~$K_n$, so 
the above results all state that the existence of some homomorphism is equivalent
to the existence of some other homomorphism. More precisely, there exists a homomorphism
of some $\Lambda(G)$ to a target $K$ if and only if there exists a homomorphism of
$G$ to some $\Gamma(K)$. In Theorem~\ref{gest}, we have $\Lambda(G) = G[K_2]$ and 
$\Gamma(K_n) = K(n,2)$. In Theorem~\ref{elsa}, $\Lambda(G) = G\times H$ and
$\Gamma(K_n) = K_n^H$ for some fixed $H$. In Theorem~\ref{albe}, $\Lambda(G) = G^{1/3}$
and $\Gamma(C_5) = K_5$. It may not yet be clear that these various $\Gamma$ all generalize
to well-defined functors. We present one more example where the presentation of
$\Gamma(K_n)$ is cryptic, and even the existence of an appropriate $\Lambda$ is not obvious.

\begin{theorem}[Gy\'arf\'as, Jensen, Stiebitz] \label{gyjest}
A graph $G$ admits an $n$-colouring in which the neighbourhood of each colour class
is an independent set if and only if $G$ admits a homomorphism to the graph $U(n)$.
\end{theorem}

It is time to introduce some relevant terminology.

\section{Pultr templates and functors}
We refer the reader to \cite{HN:Homomorphisms} for an introduction
to graph homomorphisms.
\begin{defi} 
\label{defi:pultem}
\begin{itemize} 
\item[]
\item[(i)] A {\em Pultr template} is a quadruple $\T = (P,Q,\epsilon_1,\epsilon_2)$
where $P$, $Q$ are graphs and $\epsilon_1, \epsilon_2$ homomorphisms of $P$ to $Q$ such that
$Q$ admits an automorphism $q$ with $q \circ \epsilon_1 = \epsilon_2$ and 
$q \circ \epsilon_2 = \epsilon_1$.
\item[(ii)] Given a Pultr template $\T = (P,Q,\epsilon_1,\epsilon_2)$,
the \emph{left Pultr functor}~$\Lambda_{\T}$ is the following construction:
For a graph $G$, $\Lambda_{\T}(G)$ contains one copy $P_u$ of $P$
for every vertex $u$ of $G$, and for every edge $[u,v]$ of $G$, 
$\Lambda_{\T}(G)$ contains a copy $Q_{u,v}$ of $Q$ with $\epsilon_1(P)$ identified
with $P_u$ and $\epsilon_2(P)$ identified with~$P_v$.
\item[(iii)] Given a Pultr template $\T = (P,Q,\epsilon_1,\epsilon_2)$
the \emph{central Pultr functor} $\Gamma_{\T}$ is the following construction:
For a graph $K$, the vertices of $\Gamma_{\T}(K)$ are the homomorphisms
$g : P \rightarrow K$, and the edges of $\Gamma_{\T}(K)$ are the pairs
$[g_1,g_2]$ such that there exists a homomorphism $h : Q \rightarrow K$ with
$g_1 = h \circ \epsilon_1$, $g_2 = h \circ \epsilon_2$.
\end{itemize}
\end{defi}

Note that the automorphism $q$ of $Q$ interchanging $\epsilon_1$ and $\epsilon_2$ in (i)
makes the conditions in (ii) and (iii) symmetric, so that $\Lambda_{\T}(G)$
is well defined, and $\Gamma_{\T}(K)$ is a graph rather than a digraph.
The following two examples model Theorems~\ref{gest} and \ref{elsa}.
\begin{exas}
\begin{itemize}
\item[]
\item Let $\T = (K_2,K_4,\epsilon_1,\epsilon_2)$, where $\epsilon_1, \epsilon_2$ are homomorphism
mapping $K_2$ to two non-incident edges of $K_4$. Then $\Lambda_{\T}(G)$ is the lexicographic product 
$G[K_2]$. The vertices of $\Gamma_{\T}(K_n)$ are essentially the arcs of $K_n$, that is, the couples
$(i,j)$ with $1 \leq i, j \leq n$, $i \neq j$. Its edges are the pairs $[(i,j),(i',j')]$ such that
$\{i,j\} \cap \{i',j'\} = \emptyset$. Thus  $\Gamma_{\T}(K_n)$ is the Kneser graph
$K(n,2)$ with all its vertices doubled; it is homomorphically equivalent to $K(n,2)$.
Thus Theorem~\ref{gest} states that $\Lambda_{\T}(G)$ admits a homomorphism to
$K_n$ if and only if $G$ admits a homomorphism to $\Gamma_{\T}(K_n)$.
\item Let $H$ be a fixed graph and 
$\T_H = (H \times K_1,H \times K_2,\epsilon_1,\epsilon_2)$,
where $\epsilon_1, \epsilon_2$, are the two natural injections of the independent
set $H \times K_1$ into the categorial product $H \times K_2$.
Then $\Lambda_{\T_H}(G)$ is the categorial product $G \times H$,
and $\Gamma_{\T_H}(K)$ is the exponential graph $K^H$.
Theorem~\ref{gest} states that $\Lambda_{\T_H}(G)$ 
admits a homomorphism to $K_n$ if and only if $G$ admits a homomorphism to 
$\Gamma_{\T_H}(K_n)$. 
\end{itemize}
\end{exas}
For any standard graph product $\star$ (see \cite{ImrKla:GrP})
and any graph $H$, we can similarly form a template
$\T = (H \star K_1,H \star K_2,\epsilon_1,\epsilon_2)$,
such that $\Lambda_{\T}(G) = G \star H$, and 
$\Gamma_{\T}(K)$ is an ``exponential structure'' in the sense
that $G \star H$ admits a homomorphism to $K$ if and only if 
$G$ admits a homomorphism to $\Gamma_{\T}(K)$.
Geller and Stahl \cite{GellerStahl} proved that $\chi(G[H])$ equals $\chi(G[K_{\chi(H)}])$,
and used an extension of Theorem~\ref{gest} to prove that the latter is measured by 
homomorphisms into Kneser graphs.
El-Zahar and Sauer \cite{es} used the correspondence
of Theorem~\ref{elsa} to prove that the chromatic number of a categorial 
product of $4$-chromatic graphs is $4$.  

The next examples involve Pultr functors outside the mould of graph products.
\begin{itemize}
\item Let $P_3$ denote the path with three edges and $\epsilon_1, \epsilon_2$ be the homomorphisms
mapping $K_1$ to the endpoints of $P_3$.  Let $\T_3=(K_1,P_3,\epsilon_1,\epsilon_2)$.
Then $\Lambda_{\T_3}(G)$ is obtained from $G$ by replacing each edge
by a path with three edges, that is, $G^{1/3}$.
$\Gamma_{\T_3}(K)$ is obtained from $K$ by adding edges between vertices
joined by a walk of length $3$ in $K$.
In particular, $\Gamma_{\T_3}(K)$ contains loops if and only if
$K$ contains triangles or loops, and $\Gamma_{\T_3}(C_5) = K_5$.Theorem~\ref{albe} states that $\Lambda_{\T_3}(G)$ admits a homomorphism to
$C_5$ if and only if $G$ admits a homomorphism to $\Gamma_{\T_3}(C_5))$.
\end{itemize}

Theorem~\ref{albe} is an adaptation by Hell and Ne\v{s}et\v{r}il
of a result by Maurer, Sudborough and
Welzl~\cite{MauSudWel:On-the-complexity-of-the-general}. It is an
example of a reduction among homomorphisms problems.
With similar reductions, Hell and Ne\v{s}et\v{r}il~\cite{HelNes:Dicho}
eventually proved that for any fixed non-bipartite
graph $H$, the problem of determining whether an input graph $G$ admits a homomorphism
into $H$ is NP-complete.

Theorems~\ref{gest}, \ref{elsa}, \ref{albe} are particular manifestations of the following general property:
\begin{theorem}[Pultr \cite{pultr}] \label{puthm}
For any Pultr template $\T$ and any graphs $G, K$, there exists a homomorphism
of $\Lambda_{\T}(G)$ to $K$ if and only if there exists a homomorphism
of $G$ to $\Gamma_{\T}(K)$.
\end{theorem}

\section{Adjoint functors and categories}
Two functors $\Lambda$ and $\Gamma$ are said to be respectively
{\em left} and {\em right adjoints} of each other if there is a natural correspondence
between the morphisms of $\Lambda(X)$ to $Y$ and the morphisms
of $X$ to $\Gamma(Y)$. Note that this correspondence between morphisms
is apparently a stronger statement than the existential statement of Theorem~\ref{puthm},
but this depends on the precise categorial context.
\begin{itemize}
\item[{\bf (i)}] In the usual category of graphs, it can be shown that
for the templates $\T$ modelling Theorems~\ref{gest} and~\ref{elsa},
$\Lambda_{\T}$ and $\Gamma_{\T}$ are left and right adjoints
in the sense above. However, the template $\T$ modelling Theorem~\ref{albe}
does not give rise to left and right adjoints: The number of homomorphisms
of $\Lambda_{\T}(G)$ to $K$ is not always equal to the number of homomorphisms
of $G$ to $\Gamma_{\T}(K)$
\item[{\bf (ii)}] In the ``thin'' category (or preorder) of graphs, the morphisms
between given graphs is not distinguished: there is at most
one generic morphism from one graph to another.
In this context, Theorem~\ref{puthm} states that for any Pultr template $\T$,
$\Lambda_{\T}$ and $\Gamma_{\T}$ are left and right adjoints of each other.
\item[{\bf (iii)}] In the category of multigraphs, where morphisms must specify images
of vertices and also of edges, $\Lambda_{\T}$ and $\Gamma_{\T}$
are always left and right adjoints. Furthermore Pultr~\cite{pultr} has shown that all pairs of
adjoint functors in this category are of the form $\Lambda_{\T}$ and $\Gamma_{\T}$.
\end{itemize}
Though the first context is the most commonly understood, our applications so far
and those to come are existential. Consequently we work in the thin category of graphs,
and call a pair of functors $\Lambda, \Gamma$ left and right adjoints of each other 
if the existence of a homomorphism of $\Lambda(G)$ to $K$ is equivalent to the existence 
of a homomorphism of $G$ to $\Gamma(K)$. 

Any Pultr template $\T$ gives rise to the
adjoints $\Lambda_{\T}$ and $\Gamma_{\T}$. But 
unlike the case of the category of multigraphs, other pairs of adjoint functors exist.
In particular, $\Gamma_{\T}$ is called a ``central'' rather than a ``right'' functor
because in some significant cases $\Gamma_{\T}$ itself admits a right adjoint
$\Omega_{\T}$. For instance, in the next section we interpret Theorem~\ref{gyjest}
in terms of the right adjoint of a central Pultr functor. 

We do not know which central Pultr functors admit right adjoints. 
It would be interesting to characterize all pairs of adjoint functors in the thin category
of graphs, though this objective may be out of reach. We will instead use the known 
examples to show the type of applications that the search for new adjoint functors
may yield.

\section{The right adjoint of $\Gamma_{\T_3}$}

Recall that the Pultr template $\T_3$ is
$(K_1,P_3,\epsilon_1,\epsilon_2)$, where $\epsilon_1$ and $\epsilon_2$
map $K_1$ to the endpoints of $P_3$.
Any graph $G$ is a spanning subgraph of $\Gamma_{\T_3}(G)$. 
Therefore any proper $n$-colouring of $\Gamma_{\T_3}(G)$ is 
a fortiori a proper $n$-colouring of $G$.  Let $c: G \rightarrow K_n$ 
be a proper $n$-colouring. Then $c$ is {\em not} a proper $n$-colouring
of $\Gamma_{\T_3}(G)$ if and only if $G$ contains a path of 
length three whose end vertices are identically coloured. 
The middle points of this path are then adjacent neighbours of a colour class.
Therefore, the proper $n$-colourings of $\Gamma_{\T_3}(G)$ are precisely 
the proper $n$-colourings of $G$ such that the neighbourhood of each colour class is an independent set.
Thus Theorem~\ref{gyjest} states that there exists a homomorphism of $\Gamma_{\T_3}(G)$
to $K_n$ if and only if there exists a homomorphism of $G$ to some $U(n)$. 

\begin{defi} \rm
For a graph $H$, let $\Omega_{\T_3}(H)$ be the graph constructed as follows.
The vertices of $\Omega_{\T_3}(H)$ are the couples $(u,U)$ such that
$u \in V(H)$ and $U \subseteq N_H(u)$, the neighbourhood of~$u$ in~$H$. Two couples $(u,U)$, $(v,V)$ are joined
by an edge of $\Omega_{\T_3}(H)$ if $u \in V$, $v \in U$, and every vertex
in $U$ is adjacent to every vertex in $V$.
\end{defi}

\begin{theorem}[\cite{tarmul}] \label{rap3}
For any graphs $G$ and $H$, there exists a homomorphism of $\Gamma_{\T_3}(G)$
to $H$ if and only if there exists a homomorphism of~$G$ to~$\Omega_{\T_3}(H)$.
\end{theorem}

With $H = K_n$ and $U(n) = \Omega_{\T_3}(K_n)$, this is the statement of Theorem~\ref{gyjest}.
The purpose of Gy\'arf\'as, Jensen and Stiebitz~\cite{GyaJenSti:On-graphs-with} was to answer affirmatively a question of Harvey and Murty
by showing that there exists, for every $n$, a $n$-chromatic graph with ``strongly independent colour classes'',
that is, a $n$-chromatic graph $G_n$ such that $\Gamma_{\T_3}(G_n)$ is $n$-chromatic.
By Theorem~\ref{gyjest}, such a graph exists if and only if $G_n = \Omega_{\T_3}(K_n)$ has this property.
They prove that indeed $\chi(\Omega_{\T_3}(K_n)) = n$.

The purpose in \cite{tarmul} was to find multiplicative graphs. 
\begin{defi}
A graph $K$ is {\em multiplicative} if whenever a product
$G \times H$ admits a homomorphism to $K$, one of the factors
$G$ or $H$ admits a homomorphism to $K$. 
\end{defi}
For a long time, only $K_2$ and the odd cycles were known to be multiplicative.
Adjoint functors help to find new multiplicative graphs from known ones:
For any Pultr template $\T$, we have 
$\Gamma_{\T}(G\times H) \simeq \Gamma_{\T}(G)
\times \Gamma_{\T}(H)$. Using this property, it is not hard to show
that if $\Gamma_{\T}$ admits a right adjoint~$\Omega_{\T}$,
then for any multiplicative graph $K$, $\Omega_{\T}(K)$ is multiplicative.
In the case of~$\T_3$, we have $\Gamma_{\T_3}(\Omega_{\T_3}(G))$
homomorphically equivalent to $G$ for any graph $G$, and this allows to prove the following.
\begin{theorem}[\cite{tarmul}]
For any graph $K$, $K$ is multiplicative if and only if $\Omega_{\T_3}(K)$
is multiplicative.
\end{theorem}
For relatively prime positive integers $m$, $n$ such that $2m \le n$, the \emph{circular complete graph}
$K_{n/m}$ is the graph whose vertices are the elements of the
cyclic group~$\mbox{\Bb Z}_n$, where $u$ and $v$ are joined by an edge
if $u-v \in \{m, \ldots, n-m\}$. Note that $K_{(2m+1)/m}$ is the odd cycle $C_{2m+1}$,
and for $n/m < 3$, $\Gamma_{\T_3}(K_{n/m}) \simeq K_{n/(3m-n)}$.
It can be shown that for $n/m < 12/5$, $\Omega_{\T_3}(K_{n/(3m-n)})$ is homomorphically
equivalent to $K_{n/m}$.
Using these results, it was possible to show that the circular complete graphs $K_{n/m}$
with $n/m < 4$ are all multiplicative. For a while, it looked like the same method
would yield many new discoveries of multiplicative graphs. None have yet been found, 
but the results of Hajiabolhassan and Taherkhani, which we present next,
have exhibited more links between similar functors and circular complete graphs.

\section{Odd powers and roots}
\label{sec:odd}

In this section we present generalizations of the functors
$\Lambda_{\T_3}, \Gamma_{\T_3}$ and $\Omega_{\T_3}$ studied by
Hajiabolhassan and Taherkhani \cite{HaT:Graph-powers}. 
For an integer $m$, let $P_m$ denote the path with $m$~edges.
For odd~$m$, let $\T_m =(K_1,P_m,\epsilon_1,\epsilon_2)$, where $\epsilon_1, \epsilon_2$
are the homomorphisms mapping~$K_1$ to the endpoints of~$P_m$. 
Then $\Lambda_{\T_m}(G)$ is the graph obtained from~$G$ by replacing each edge by
a copy of~$P_m$, that is, the $m$-subdivision~$G^{1/m}$ of~$G$. 
$\Gamma_{\T_m}(H)$ is the ``$m$-th power of~$H$'', obtained
from~$H$ by adding edges between pairs of vertices connected by a walk
of length~$m$. (In particular, $\Lambda_{\T_1}(G) = \Gamma_{\T_1}(G) = G$.)
We now describe a right adjoint of~$\Gamma_{\T_m}$: 
\begin{defi} \rm
For an odd integer $m = 2k+1$ and a graph $H$, 
let $\Omega_{\T_m}(H)$ be the graph constructed as follows.
The vertices of $\Omega_{\T_m}(H)$ are the $(k+1)$-tuples 
$(u,U_1, \ldots, U_k)$ such that
$u \in V(H)$, $U_1 \subseteq N_H(u)$, $U_i \subseteq V(H)$ and $U_i$ is completely joined to $U_{i-1}$
for $i = 2, \ldots, k$. Two $k$-tuples $(u,U_1, \ldots, U_k)$, 
$(v,V_1, \ldots, V_k)$ are joined
by an edge of $\Omega_{\T_m}(G)$ if $u \in V_1$, $v \in U_1$, 
$U_{i-1} \subseteq V_i$ and $V_{i-1} \subseteq U_i$ for $i = 2, \ldots, k$,
and $U_k$ is completely joined to $V_k$. (Here, $N_H(u)$~is the set
of all vertices of~$H$ adjacent to~$u$, and two sets of vertices
are called {\em completely joined} if every vertex in one is adjacent
to every vertex in the other.)
\end{defi}

Theorem~\ref{rap3} generalizes as follows.
\begin{theorem}[\cite{HaT:Graph-powers}]
\label{thm:omega}
For two graphs $G$ and $H$, there exists a homomorphism of $\Gamma_{\T_m}(G)$ to $H$
if and only there exists a homomorphism of $G$ to $\Omega_{\T_m}(H)$.
\end{theorem}

For odd $s$ and $r$, define $P^s_r(G) = \Gamma_{\T_s} ( \Lambda_{\T_r} (G))$
and $R^r_s(H) = \Gamma_{\T_r} ( \Omega_{\T_s} (H))$. 
There exists a homomorphism of $P^s_r(G) = \Gamma_{\T_s} ( \Lambda_{\T_r} (G))$ 
to $H$ if and only if
there exists a homomorphism of $\Lambda_{\T_r} (G)$ to $\Omega_{\T_s}(H)$, 
that is, if and only if there exists a homomorphism
of $G$ to $\Gamma_{\T_r} ( \Omega_{\T_s}(H) ) = R^r_s(H)$.
Thus $P^s_r$ and $R^r_s$ are right and left adjoint of each other, though
they are not necessarily left, central or right functors 
associated to Pultr templates. These are ``ordered'' as follows.
\begin{theorem}[\cite{HaT:Graph-powers}]
\label{thm:5.3}
Let $s, r, s', r'$ be odd integers such that $\frac{s}{r} \leq \frac{s'}{r'}$.
Then for any graph $G$, $P^s_r(G)$ admits a homomorphism to
$P^{s'}_{r'}(G)$ and $R^{r'}_{s'}(G)$ admits a homomorphism to
$R^{r}_{s}(G)$.
\end{theorem}
The {\em circular chromatic number} $\chi_{\rm c}(G)$ of a graph $G$ is
the minimum value $n/m$ such that $G$ admits a homomorphism to the circular complete
graph $K_{n/m}$. Note that for odd $s = 2i+1$, $\Omega_{\T_s}(K_3)$ is homomorphically equivalent
to the $3s$-cycle $K_{(6i+3)/(3i+1)}$, and for $r = 2j+1$, $R^r_s(K_3) = \Gamma_{\T_r}(\Omega_{\T_s}(K_3))$
is homomorphically equivalent to $K_{(6i+3)/(3i+1-j)}$. Using the fact that
$P^s_r$ and $R^r_s$ are right and left adjoint of each other, we get the following.
\begin{theorem}[\cite{HaT:Graph-powers}]
For a graph $G$, $\chi_{\rm c}(G)$ is the supremum of the values
$(6i+3)/(3i+1-j)$ such that $P^{2i+1}_{2j+1}(G)$ is $3$-colourable.
\end{theorem}

\section{Oriented paths as templates}
\label{sec:opa}

In this and the following sections we change the setting and consider
digraphs rather than undirected (symmetric) graphs. To this end,
we need to modify slightly the definition of a Pultr template.

\begin{defi}
In the setting of digraphs, a \emph{Pultr template} is a quadruple
$\T=(P,Q,\epsilon_1,\epsilon_2)$, where $P,Q$~are digraphs and
$\epsilon_1,\epsilon_2$~are homomorphisms of~$P$ to~$Q$.
\end{defi}

Thus we no longer require the existence of a special automorphism~$q$
of~$Q$, whose purpose was to ensure that $\Gamma_\T(G)$~would be
an undirected graph for any undirected graph~$G$.

To define the left Pultr functor and the central Pultr functor
corresponding to a Pultr template~$\T$, simply replace the word
``edge(s)'' with ``arc(s)'' in Definition~\ref{defi:pultem} (ii), (iii).

In the rest of this section, we present an oriented analogue of the
construction of~$\Omega_\T$ and of Theorem~\ref{thm:omega}.

Let $Q$ be an orientation of a path and consider the Pultr template
$\T=(K_1,Q,\epsilon_1,\epsilon_2)$, where $\epsilon_1$, $\epsilon_2$
are the homomorphisms mapping $K_1$ to the end-points of~$Q$.
Similarly to the situation in Section~\ref{sec:odd}, $\Lambda_\T(G)$~is
the digraph obtained from~$G$ by replacing each arc with a copy
of the path~$Q$. $\Gamma_\T(H)$~is the digraph on the same vertex set as~$H$
in which there is an arc from~$u$ to~$v$ if and only if there exists
an oriented walk from~$u$ to~$v$ in~$H$ whose steps are oriented
according to the orientations of the arcs of~$Q$. The right adjoint
of~$\Gamma_\T$ is as follows:

\begin{defi}
Suppose $Q$ has vertices $0,1,\dotsc,m$. For a digraph~$H$, let
$\Omega_\T(H)$ be the following digraph: The vertices of~$\Omega_\T(H)$
are all the $({m+1})$-tuples $(u,U_1,\dotsc,U_{m})$ such that $u\in V(H)$
and $U_i\subseteq V(H)$ for $i=1,2,\dotsc,m$, with $u\Rrightarrow U_m$.  There is an arc
in~$\Omega_\T(H)$ from $(u,U_1,\dotsc,U_{m})$ to $(v,V_1,\dotsc,V_{m})$
if and only if
\begin{itemize}
\item[(1a)] $u\in V_1$ if $0\to 1$ in $Q$,
\item[(1b)] $v\in U_1$ if $1\to 0$ in $Q$; and
\item[(2)] for each $i=1,\dotsc,m-1$:
\begin{itemize}
\item[(a)] $U_i\subseteq V_{i+1}$ if $i\to i+1$ in $Q$,
\item[(b)] $V_i\subseteq U_{i+1}$ if $i+1\to i$ in $Q$.
\end{itemize}
\end{itemize}
(The notation $a\to b$ means that there is an arc from~$a$ to~$b$ and
$u \Rrightarrow V$ means that there is an arc from~$u$ to every element of~$V$.)
\end{defi}

\begin{theorem}[\cite{FonTar:Right-adjoints}]
\label{thm:omega-path}
For any two digraphs $G$ and $H$, there exists a homomorphism
of~$\Gamma_\T(G)$ to~$H$ if and only if there exists a homomorphism
of~$G$ to~$\Omega_\T(H)$.
\end{theorem}

With oriented paths, the structure and ordering of functors analogous
to~$P_r^s$ and~$R_s^r$ (see Section~\ref{sec:odd}) gets much more
complex. Possible applications are currently unknown (see also
Section~\ref{sec:op}).

In the next section, we consider Pultr functors for which, on
the other hand, an abundance of applications can be found in the
literature.

\section{Shift graphs}
\label{sec:shift}

\begin{defi}
Let $H$ be a digraph. The \emph{arc graph} of~$H$ is the
digraph~$\delta(H)$ whose vertices are the arcs of~$H$ and
$(u,v)\to(x,y)$ in~$\delta(H)$ if and only if $v=x$.
\end{defi}

Observe that $\delta$ is the central Pultr functor given by the
template $(\vec P_1,\vec P_2,\epsilon_1,\epsilon_2)$, where
\begin{align*}
&\begin{aligned}
\vec P_1 &= 0\to 1, \\
\vec P_2 &= 0\to 1\to 2,
\end{aligned}\\
&\epsilon_1: i\mapsto i, \qquad \epsilon_2: i\mapsto i+1.
\end{align*}

A convenient fact about arc graphs is that we know good bounds on
the chromatic number of~$\delta(H)$ in terms of the chromatic number
of~$H$ (see~\cite{HarEnt:Arc-colorings,PolRod:Arc}). Given a proper
$k$-colouring of~$\delta(H)$, we can construct
a proper $2^k$-colouring of~$H$: let the colour of a vertex~$u$
of~$H$ be the set of all colours used on the outgoing arcs from~$u$
in the proper $k$-colouring of~$\delta(H)$. This shows that
$\chi(\delta(H))\ge\log\chi(H)$. (In fact, we have $\chi(\delta(H)) = \Theta(\log\chi(H))$.)

This fact was used by Poljak and R\"odl~\cite{PolRod:Arc} to discuss
the possible boundedness of the ``Poljak-R\"odl function'': Define
$f : \mathbb{N} \to \mathbb{N}$ by letting $f(n)$ be the minimum
possible chromatic number $\chi(G\times H)$ of a product of 
$n$-chromatic digraphs $G$ and $H$. (Here, the chromatic number of a digraph is
defined to be the chromatic number of its symmetrization.)
It is not known whether $f$ is bounded or unbounded, but
Poljak and R\"odl were able to limit the possible upper bounds on $f$.
One element of the proof is the bound $\chi(\delta(G))\ge\log\chi(G)$
and the other is the fact that $\delta$, like all central Pultr functor,
commutes with the categorial product, that is, the identity
$\delta(G\times H) \simeq \delta(G) \times \delta(H)$.
The best result in this direction is the following:
\begin{theorem}[\cite{Pol:cdia}]
Either $f(n) \leq 3$ for all $n$, or $\lim_{n \to \infty}f(n) = \infty$.
\end{theorem}
Let $g$ be the undirected analogue of $f$. Hedetniemi's conjecture states that
$g(n) = n$ for all $n$. Using the results on directed graphs,
it is possible to prove that $g$ is either bounded above by $9$ 
or unbounded (see~\cite{Pol:cdia}). In fact, $g$ is unbounded
if and only if $f$ is unbounded (see~\cite{TarWeh:prf}).

Properties of the arc graph construction are also used in the analysis of
the ``shift graphs'' that are the folklore examples of graphs
with large odd girth and large chromatic number.

\begin{defi}
Let $n$, $k$ be positive integers, $k\ge 2$. The \emph{directed
shift graph} $R(n,k)$ is the digraph with vertex set
$V(R(n,k))=\{(u_1,\dotsc,u_k)\colon 1\le u_1< u_2 < \dotsb <u_k\le
n\}$ where $(u_1,\dotsc,u_k)\to (v_1,\dotsc,v_k)$ if and only if
$u_2=v_1$, $u_3=v_2$, \dots, $u_k=v_{k-1}$.

The \emph{undirected shift graph} $R'(n,k)$ is the symmetrization
of~$R(n,k)$.
\end{defi}
\begin{theorem}[\cite{NesRod:Type-theory}]
Let $c,k\ge 2$ and put $n=2^{2^{\cdots 2^c}}$, where the tower of
powers has height~$k$. Then the undirected shift graph $R'(n,k)$
has chromatic number at least~$c$ and odd girth at least $2k+1$.
\end{theorem}
We will show that both properties of high chromatic number and high odd girth
are related to properties of adjoint functors.
First note that $R(n,k)=\Gamma_{\T_k}(\vec T_n)$ for the Pultr
template $\T_k=(\vec P_k,\vec P_{k+1},\epsilon_1,\epsilon_2)$ with
\begin{align*}
&\begin{aligned}
\vec P_k &= 0\to 1 \to \dotsb \to k-1, \\
\vec P_{k+1} &= 0\to 1 \to \dotsb \to k,
\end{aligned}\\
&\epsilon_1: i\mapsto i, \qquad \epsilon_2: i\mapsto i+1.
\end{align*}
At the same time, 
$R(n,2)=\delta(\vec T_n)$ and for $k\ge 3$ we have
$R(n,k)\simeq \delta(R(n,k-1))$. Hence we can get all shift graphs by
iterating the arc graph functor, starting with a transitive tournament.
The bound $\chi(\delta(G))\ge\log\chi(G)$ directly implies that 
$\chi(R(n,k))\ge\log^{k-1}n$,
where $\log^{k-1}$ means the binary logarithm iterated $k-1$ times.

Next we are going to show that the odd girth of the undirected shift
graphs is large, namely the odd girth of $R'(n,k)$ is at least
$2k+1$. Suppose that some odd cycle $C$ admits a homomorphism
to $R'(n,k)$. Then there exists an orientation $\vec{C}$ of $C$
which admits a homomorphism to $R(n,k) \simeq \delta(R(n,{k-1}))$.
Therefore there exists a homomorphism of~$\delta_L(\vec{C})$
to~$R(n,{k-1})$, where $\delta_L$~is the left adjoint of~$\delta$.
By construction,  $\delta_L(\vec{C})$ contains an arc ${u_0 \to u_1}$
for every vertex $u$ of $\vec{C}$, with $u_1$ identified to
$v_0$ for every arc $u \to v$ of $\vec{C}$. Thus the number
of vertices of $\delta_L(\vec{C})$ is the same as that of $\vec{C}$.
Also, $\delta_L(\vec{C})$ is not bipartite since
a homomorphism  of $\delta_L(\vec{C})$
to $K_2$ would correspond to a homomorphism
of $\vec{C}$ to $\delta(K_2) \simeq K_2$, which is impossible since $C$
is an odd cycle. Therefore, $\delta_L(\vec{C})$ contains an odd cycle.
However, $R(n,k)$ has no directed cycles (since the projection
on the first coordinate is a homomorphism to a transitive tournament)
hence $\vec{C}$ has at least one source $u$ and one sink $v$.
The two vertices $u_0$ and $v_1$ are then respectively a source
and a sink in $\delta_L(\vec{C})$. Hence the odd girth of (the symmetrization of)
$\delta_L(\vec{C})$ is smaller than that of $C$. Since the odd girth
of $K_n = R(n,1)$ is $3$, this implies that the odd girth of
$R(n,k)$ is at least $2k + 1$. (In fact, the odd girth $R(n,k)$ 
of is exactly $2k + 1$, unless $n < 2k+1$.) 

\section{Pultr functors and tree duality}

\begin{defi}
A set $\F$ of digraphs is a \emph{complete set of obstructions} for
a digraph~$H$ if for any digraph~$G$ there exists \emph{no}
homomorphism of~$G$ to~$H$ if and only if there exists a homomorphism
of \emph{some} $F\in\F$ to~$G$. We also say that $(\F,H)$ is a
\emph{homomorphism duality}.
\end{defi}

If $H$ admits a finite complete set of obstructions, then we say
that $H$~\emph{has finite duality}; in this case,
by~\cite{Kom:Phd,NesTar:Dual}, $H$~admits a finite complete set of
obstructions all of whose elements are oriented trees.  Conversely,
every finite set~$\F$ of oriented trees is a complete set of
obstructions for some digraph~$H$. If $H$~admits a (not necessarily
finite) complete set of obstructions all of whose elements are
oriented trees, we say that $H$~\emph{has tree duality}.

In~\cite{FonTar:Adjoint} we proved that if $H$~has tree duality,
then so does its arc graph~$\delta(H)$. Furthermore we gave an
explicit description of a complete set of tree obstructions
for~$\delta(H)$, provided we are given a complete set of tree
obstructions for~$H$.

\begin{defi}
Let $T$ be a tree. For every vertex $u$ of~$T$, let $F(u)$ be a
tree that admits a homomorphism to~$\vec P_1$; fix such a homomorphism
$\phi_u:F(u)\to P_1$, so that for any arc $(x,y)$ of~$F(u)$ we have
$\phi_u(x)=0$, $\phi_u(y)=1$. For each arc~$e$ of~$T$, incident
with~$u$, fix a vertex $v(u,e)$ of~$F(u)$ in such a way that if
$e$~is outgoing from~$u$, then $\phi_u(v(u,e))=1$, and
if $e$~is incoming to~$u$, then $\phi_u(v(u,e))=0$.
Construct a tree~$S$ by taking all the trees~$F(u)$ for all the
vertices~$u$ of~$T$, and by identifying the vertex $v(u,e)$ with
$v(u',e)$ for every arc $e=(u,u')$ of~$T$. Any tree~$S$ constructed
from~$T$ by the above procedure is called a \emph{sproink} of~$T$.
\end{defi}

\begin{theorem}[\cite{FonTar:Adjoint}]
\label{thm:sproink}
If $\F$ is a complete set of tree obstructions for some digraph~$H$,
then the set of all sproinks of all the trees in~$\F$ is a complete
set of tree obstructions for its arc graph~$\delta(H)$.
\end{theorem}

\begin{exa}
Let $\vec P_k$ be the directed path with $k$~arcs (that is, the
path $0\to1\to2\to\dotsb\to k$) and let $\vec T_k$ be the transitive
tournament on $k$~vertices.  By~\cite{NesPul:SubFac}, $\{\vec
P_k\}$~is a complete set of tree obstructions for~$\vec T_k$. To
get a complete set of obstructions for~$\delta(\vec T_k)$, we can
take just the \emph{minimal} sproinks of~$\vec P_k$ (minimal with
respect to the ordering by existence of homomorphisms). In the
minimal sproinks, $F(0)$ and $F(k)$ will each be the one-vertex graph~$K_1$.
All the other $F(u)$'s will be alternating paths (``zigzags'') that
will connect at their end-points. So all the minimal sproinks
of~$\vec P_k$ for $k\ge3$ can be described by the regular expression
\[ {\uparrow} ({\uparrow}({\downarrow}{\uparrow})^{\ast})^{k-3}{\uparrow} .\]
\end{exa}

Note that if $(\F,H)$ and $(\F',H')$ are dualities, then so is
$(\F\cup\F',H\times H')$. Hence starting with graphs with finite
duality, whose structure is rather well understood, and taking iterated
arc graphs and products yields a fairly large class of graphs with
tree duality as well as a complete set of obstructions for each of
them. In particular, the knowledge of a complete set of obstructions
of the directed shift graphs was used in~\cite{Tar:dbl} to prove 
the density of the lattices
$K_n^{\mathcal{D}}$ of directed graph powers of $K_n$
(under homomorphic equivalence).

There is nothing special about the arc graph, however. In fact, all
central Pultr functors and all their right adjoints preserve tree
duality:
\begin{theorem}[\cite{FonTar:Adjoint,FonTar:Right-adjoints}]
Let $\T$ be a Pultr template and let $H$ be a digraph with tree
duality. Then $\Gamma_\T$~has tree duality. Furthermore, if there
exists a digraph $\Omega_\T(H)$ such that, for any digraph~$G$,
$G\to\Omega_\T(H)$ iff $\Gamma_\T(G)\to H$, then $\Omega_\T(H)$~has
tree duality.
\end{theorem}

The tree obstructions for~$\Omega_\T(H)$ have a neat description
using the left adjoint~$\Lambda_\T$. On the other hand, an explicit
description of the tree obstructions for~$\Gamma_\T(H)$ for a general
Pultr template~$\T$ is currently unknown. The knowledge of the
obstructions in some special cases can have interesting applications,
as we show next.

\section{Circular Gallai--Roy theorem}

The following well-known theorem is usually credited to
Gallai~\cite{Gal:On-directed} and Roy~\cite{Roy:Nombre}, even though
it had independently been proved earlier by Vitaver~\cite{Vit:Determination} and
Hasse~\cite{Has:Zur-algebraischen}.

\begin{theorem}[Vitaver, Hasse, Roy, Gallai]
\label{thm:vhrg}
A graph $G$ is $k$-colourable if and only if it admits an
orientation~$\vec G$ such that there is no homomorphism of~$\vec
P_k$ to~$\vec G$.
\end{theorem}

As we have already mentioned, $\{\vec P_k\}$~is a complete set of
obstructions for~$\vec T_k$, the transitive tournament on $k$~vertices.
That is, for any digraph~$\vec G$ we have $\vec P_k\notto\vec G$
if and only if $\vec G\to\vec T_k$. Observe that $G$~is $k$-colourable
if and only if it admits an orientation~$\vec G$ such that $\vec
G\to\vec T_k$ and you get Theorem~\ref{thm:vhrg}.

We are now interested in finding an analogous condition for
circular colourability.
Recall that for relatively prime integers $m, n$ such that $2m \le
n$, the circular complete graph $K_{n/m}$ is the graph whose vertices
are the elements of the cyclic group $\mbox{\Bb Z}_n$, where $u$
and $v$ are joined by an edge if $u-v \in \{m, \ldots, n-m\}$; the
circular chromatic number $\chi_{\rm c}(G)$ of a graph $G$ is the
minimum value $n/m$ such that $G$ admits a homomorphism to the
circular complete graph $K_{n/m}$.

\begin{defi}
Let $m\ge 1$ be an integer. The \emph{$m$-th interleaved adjoint}
of a digraph~$H$ is the digraph~$\iota_m(H)$ whose vertices are all
the $m$-tuples of vertices of~$H$, and $(u_1,\dotsc,u_m)\to
(v_1,\dotsc,v_m)$ in~$\iota_m(H)$ if $u_i\to v_i$ in~$H$ for all
$i=1,\dotsc,m$ and $v_i\to u_{i+1}$ in~$H$ for all $i=1,\dotsc,m-1$.
\end{defi}

It turns out that $\iota_m$~is a central
Pultr functor with the template $(P_m,Q_m,\epsilon_{m,1},\epsilon_{m,2})$,
where $P_m$ has vertices $1,2,\dotsc,m$ and no arcs; $Q_m$~has
vertices $1_1,1_2,2_1,2_2,\dotsc,m_1,m_2$ and arcs $u_1\to u_2$ for
$u=1,2,\dotsc,m$ and $u_2\to (u+1)_1$ for $u=1,2,\dotsc,m-1$, and
$\epsilon_{m,1}(u)=u_1$, $\epsilon_{m,2}(u)=u_2$ for all $u=1,2,\dotsc,m$.

Let $\lambda_m$ be the left adjoint of $\iota_m$. Hence in particular
$\lambda_m(G)\to\vec T_n$ if and only if $G\to\iota_m(\vec T_n)$.
By homomorphism duality, $\lambda_m(G)\to\vec T_n$ if and only if
$\vec P_n\notto\lambda_m(G)$. Combining these two equivalences and
using the explicit description of~$\lambda_m(G)$ (given by
Definition~\ref{defi:pultem}(ii)) we can get the following.

\begin{prop}[\cite{FonNesTar:Interlacing}]
Let $\P_{n,m-1}$ be the family of oriented paths obtained from the
directed path~$\vec P_n$ by reversing at most $m-1$ arcs.  For a
digraph~$G$, there exists a homomorphism of~$G$ to~$\iota_m(\vec
T_n)$ if and only if there exists no homomorphism to~$G$ from any
path in~$\P_{n,m-1}$.
\end{prop}

A surprising connection between circular colourings and interleaved
adjoints has been discovered by Yeh and
Zhu~\cite{YehZhu:Resource-sharing-system}.

\begin{theorem}[\cite{YehZhu:Resource-sharing-system}]
For integers $m,n$ such that $n\ge 2m$, there exist homomorphisms
both ways between $K_{n/m}$ and $B(n,m)$, the symmetrization
of~$\iota_m(\vec T_n)$.
\end{theorem}

Thus, for an undirected graph~$G$ there exists a homomorphism of~$G$
to~$K_{n/m}$ if and only if there exists an orientation~$\vec G$
of~$G$ that admits a homomorphism to~$\iota_m(\vec T_n)$.
Hence we get the sought analogue of Theorem~\ref{thm:vhrg}.

\begin{theorem}
Let $\P_{n,m-1}$ be the family of oriented paths obtained from the
directed path~$\vec P_n$ by reversing at most $m-1$ arcs.
A graph~$G$ has circular chromatic number at most $n/m$ if and only
if it admits an orientation~$\vec G$ such that there exists no
homomorphism to~$\vec G$ from any path in~$\P_{n,m-1}$.
\end{theorem}

\section{Open problems}
\label{sec:op}

Finally, we present several open problems hoping to stimulate
interest in the topic. The problem we currently find quite intriguing
and at the same time within reach is this:

\begin{problem}
For which Pultr templates~$\T$ does the central Pultr functor~$\Gamma_\T$
admit a right adjoint? This problem is open in both the directed
and the undirected case, and has a different flavour in each.
\end{problem}

The proof of Theorem~\ref{thm:5.3} makes use of the fact that
$\Gamma_{\T_m}(\Lambda_{\T_m}(G))$ is homomorphically equivalent
to~$G$ for any graph~$G$, for the path templates~$\T_m$ of
Section~\ref{sec:odd} with odd~$m$. In fact, any $G$ admits a
homomorphism to $\Lambda(\Gamma(G))$ for any pair of adjoint functors
$\Lambda, \Gamma$. So the important property of the path templates
is that $\Gamma_{\T_m}(\Lambda_{\T_m}(G))\to G$ for any graph~$G$.
This leads to the following question:

\begin{problem}
For what Pultr templates~$\T$ does $\Gamma_\T(\Lambda_\T(G))$ admit
a homomorphism to~$G$ for any~$G$?
\end{problem}

By Theorem~\ref{thm:omega-path}, the central Pultr functor~$\Gamma_\T$
admits a right adjoint~$\Omega_\T$ for path templates~$\T$ also in
the setting of digraphs. Thus we may consider directed analogues
of the functors $P_r^s$ and~$R_s^r$ of Section~\ref{sec:opa}.
However, the ordering of the path templates is no longer linear,
nor is the ordering of the corresponding $P_r^s$'s and~$R_s^r$'s.

\begin{problem}
Let $\T_1,\T_2,\T_3,\T_4$ be Pultr templates for digraphs such that
each $P$ is $K_1$, each $Q$~is an oriented path and each
$\epsilon_1,\epsilon_2$ map~$K_1$ to the end points of~$Q$. Then
by Theorem~\ref{thm:omega-path} there exists a right
adjoint~$\Omega_{\T_i}$ for each~$i$. Define
$P_i^j(G)=\Gamma_{\T_j}(\Lambda_{\T_i}(G))$; put
$R_j^i(H)=\Gamma_{\T_i}(\Omega_{\T_j}(H))$ for $i,j\in\{1,2,3,4\}$.
Under what conditions on the templates do we get an analogue of
Theorem~\ref{thm:5.3}?
\end{problem}

The following problem is motivated by the chromatic properties of
arc graphs, see Section~\ref{sec:shift}.

\begin{problem}
Characterize Pultr templates~$\T$ for which there exists an unbounded
function $c:\Nn\to\Nn$ such that $\chi(\Gamma_\T(H))\ge c(\chi(H))$
for any digraph~$H$.
\end{problem}

Finally, we would like to see a construction similar to the sproinks
of Theorem~\ref{thm:sproink}, for arbitrary Pultr templates.

\begin{problem}
Describe a complete set of tree obstructions for~$\Gamma_\T(H)$ in
terms of the tree obstructions for~$H$, for a general Pultr
template~$\T$.
\end{problem}


\providecommand{\bysame}{\leavevmode\hbox to3em{\hrulefill}\thinspace}
\providecommand{\MR}{\relax\ifhmode\unskip\space\fi MR }
\providecommand{\MRhref}[2]{%
  \href{http://www.ams.org/mathscinet-getitem?mr=#1}{#2}
}
\providecommand{\href}[2]{#2}


\begin{thebibliography}{10}

\bibitem{EZaSau:Chro}
Mohamed El-Zahar and Norbert Sauer, \emph{The chromatic number of the product
  of two 4-chromatic graphs is 4}, Combinatorica \textbf{5} (1985), no.~2,
  121--126.

\bibitem{FonNesTar:Interlacing}
Jan Foniok, Jaroslav Ne{\v s}et{\v r}il, and Claude Tardif, \emph{Interleaved
  adjoints of directed graphs}, European J. Combin. \textbf{32} (2011), no.~7,
  1018--1024.

\bibitem{FonTar:Right-adjoints}
Jan Foniok and Claude Tardif, \emph{Digraph functors which admit both left and
  right adjoints}, In preparation.

\bibitem{FonTar:Adjoint}
\bysame, \emph{Adjoint functors and tree duality}, Discrete Math. Theor.
  Comput. Sci. \textbf{11} (2009), no.~2, 97--110.

\bibitem{Gal:On-directed}
Tibor Gallai, \emph{On directed paths and circuits}, Theory of Graphs (Proc.
  Colloq., Tihany, 1966) (New York), Academic Press, 1968, pp.~115--118.

\bibitem{GellerStahl}
Dennis Geller and Saul Stahl, \emph{The chromatic number and other functions of
  the lexicographic product}, J. Combin. Theory Ser. B \textbf{19} (1975),
  no.~1, 87--95.

\bibitem{GyaJenSti:On-graphs-with}
Andr{\'a}s Gy{\'a}rf{\'a}s, Tommy~R. Jensen, and Michael Stiebitz, \emph{On
  graphs with strongly independent color-classes}, J. Graph Theory \textbf{46}
  (2004), no.~1, 1--14.

\bibitem{HaT:Graph-powers}
Hossein Hajiabolhassan and Ali Taherkhani, \emph{Graph powers and graph
  homomorphisms}, Electron. J. Combin. \textbf{17} (2010), R17, 16 pages.

\bibitem{HarEnt:Arc-colorings}
C.~C. Harner and Roger~C. Entringer, \emph{Arc colorings of digraphs}, J.
  Combinatorial Theory Ser. B \textbf{13} (1972), 219--225.

\bibitem{Has:Zur-algebraischen}
Maria Hasse, \emph{Zur algebraischen {B}egr{\"u}ndung der {G}raphentheorie.
  {I}}, Math. Nachr. \textbf{28} (1964/1965), no.~5--6, 275--290.

\bibitem{HelNes:Dicho}
Pavol Hell and Jaroslav Ne{\v s}et{\v r}il, \emph{On the complexity of
  {$H$}-coloring}, J. Combin. Theory Ser. B \textbf{48} (1990), no.~1, 92--110.

\bibitem{HelNes:GrH}
\bysame, \emph{Graphs and homomorphisms}, Oxford Lecture Series in Mathematics
  and Its Applications, vol.~28, Oxford University Press, 2004.

\bibitem{ImrKla:GrP}
Wilfried Imrich and Sandi Klav{\v z}ar, \emph{Product graphs: Structure and
  recognition}, Wiley Series in Discrete Mathematics and Optimization,
  Wiley-Interscience, New York, 2000.

\bibitem{Kom:Phd}
Pavel Kom{\'a}rek, \emph{Good characterisations in the class of oriented
  graphs}, Ph.D. thesis, Czechoslovak Academy of Sciences, Prague, 1987, In
  Czech (Dobr{\'e} charakteristiky ve t{\v r}{\'\i}d{\v e} orientovan{\'y}ch
  graf{\r u}).

\bibitem{MauSudWel:On-the-complexity-of-the-general}
Hermann~A. Maurer, Ivan~Hal Sudborough, and Emo Welzl, \emph{On the complexity
  of the general coloring problem}, Inform. and Control \textbf{51} (1981),
  no.~2, 128--145.

\bibitem{NesPul:SubFac}
Jaroslav Ne{\v s}et{\v r}il and Ale{\v s} Pultr, \emph{On classes of relations
  and graphs determined by subobjects and factorobjects}, Discrete Math.
  \textbf{22} (1978), no.~3, 287--300.

\bibitem{NesRod:Type-theory}
Jaroslav Ne{\v s}et{\v r}il and Vojt{\v e}ch R{\"o}dl, \emph{Type theory of
  partition properties of graphs}, Recent Advances in Graph Theory (Prague)
  (M.~Fiedler, ed.), Academia, 1975, pp.~405--412.

\bibitem{NesTar:Dual}
Jaroslav Ne{\v s}et{\v r}il and Claude Tardif, \emph{Duality theorems for
  finite structures (characterising gaps and good characterisations)}, J.
  Combin. Theory Ser. B \textbf{80} (2000), no.~1, 80--97.

\bibitem{Pol:cdia}
Svatopluk Poljak, \emph{Coloring digraphs by iterated antichains}, Comment.
  Math. Univ. Carolin. \textbf{32} (1991), no.~2, 209--212.

\bibitem{PolRod:Arc}
Svatopluk Poljak and Vojt{\v e}ch R{\"o}dl, \emph{On the arc-chromatic number
  of a digraph}, J. Combin. Theory Ser. B \textbf{31} (1981), no.~2, 190--198.

\bibitem{Pul:The-right-adjoints}
Ale{\v s} Pultr, \emph{The right adjoints into the categories of relational
  systems}, Reports of the Midwest Category Seminar, IV (Berlin), Lecture Notes
  in Mathematics, vol. 137, Springer, 1970, pp.~100--113.

\bibitem{Roy:Nombre}
Bernard Roy, \emph{Nombre chromatique et plus longs chemins d'un graphe}, Revue
  fran{\c c}aise d'informatique et de recherche op{\'e}rationnelle \textbf{1}
  (1967), no.~5, 129--132.

\bibitem{Tar:Mul}
Claude Tardif, \emph{Multiplicative graphs and semi-lattice endomorphisms in
  the category of graphs}, J. Combin. Theory Ser. B \textbf{95} (2005), no.~2,
  338--345.

\bibitem{Tar:Hedetniemis}
\bysame, \emph{Hedetniemi's conjecture and dense {B}oolean lattices}, Order
  \textbf{28} (2011), no.~2, 181--191.

\bibitem{TarWeh:prf}
Claude Tardif and David Wehlau, \emph{Chromatic numbers of products of graphs:
  the directed and undirected versions of the {P}oljak--{R}{\"o}dl function},
  J. Graph Theory \textbf{51} (2006), no.~1, 33--36.

\bibitem{Vit:Determination}
L.~M. Vitaver, \emph{Determination of minimal coloring of vertices of a graph
  by means of {B}oolean powers of the incidence matrix}, Dokl. Akad. Nauk SSSR
  \textbf{147} (1962), 758--759, In Russian.

\bibitem{YehZhu:Resource-sharing-system}
Hong-Gwa Yeh and Xuding Zhu, \emph{Resource-sharing system scheduling and
  circular chromatic number}, Theoret. Comput. Sci. \textbf{332} (2005),
  no.~1--3, 447--460.

\end{thebibliography}
\end{document}